\documentclass[a4paper, 11pt]{article}
\usepackage{anysize}
\marginsize{3,5cm}{2,5cm}{2,5cm}{2,5cm}
\usepackage{amssymb}
\usepackage{amsmath,amsbsy,amssymb,amscd}
\usepackage{t1enc}
\usepackage[cp1250]{inputenc}
\usepackage[british]{babel}
\usepackage[all]{xy}
\usepackage{color}
\usepackage{amsfonts}
\usepackage{latexsym}
\usepackage{amsthm}
\usepackage{mathrsfs}

\usepackage{mathrsfs} 

\DeclareMathAlphabet{\mathpzc}{OT1}{pzc}{L}{it} 

\newtheorem{definition}{Definition}[section]

\newtheorem{proposition}[definition]{Proposition}
\newtheorem{theorem}[definition]{Theorem}
\newtheorem{corollary}[definition]{Corollary}
\newtheorem{remark}[definition]{Remark}
\newtheorem{lemma}[definition]{Lemma}

\def\C{\mathbb{C}}
\def\geq{\geqslant}
\def\leq{\leqslant}
\def\R{\mathbb{R}}
\def\T{\mathbb{T}}

\def\cl{\mathrm{cl\,}}

\def\Z{\mathbb{Z}}
\def\N{\mathbb{N}}

\def\sgn{\mathrm{sgn}\,}

\def\Q{\mathbb{Q}}
\newcommand{\bea}{\begin{eqnarray}}
  \newcommand{\eea}{\end{eqnarray}}
  \newcommand{\beab}{\begin{eqnarray*}}
  \newcommand{\eeab}{\end{eqnarray*}}

  \newcommand{\be}{\begin{equation}}
  \newcommand{\ee}{\end{equation}}
\title{Ratner's property for special flows over irrational rotations under functions of bounded variation}
\author{Adam Kanigowski}
\begin{document}
\baselineskip=14pt \maketitle
\begin{abstract} We consider special flows over the rotation by an irrational $\alpha$ under the roof functions of bounded variation without continuous, singular part in the Lebesgue decomposition and the sum of jumps $\neq 0$. We show that all such flows are weakly mixing. Under the additional assumption that $\alpha$ has bounded partial quotients, we study weak Ratner's property. We establish this property whenever an additional condition (stable under sufficiently small perturbations) on the set of jumps is satisfied. While it is classical that the flows under consideration are not mixing, one more condition on the set of jumps turns out to be sufficient to obtain the absence of partial rigidity, hence mild mixing of such flows.
\end{abstract}

\section{Introduction} \indent In the paper, we deal with measurable, measure-preserving flows on standard probability spaces. Hence, we are given a standard probability space $(X,\mathscr{B},\mu)$ and an $\R$-representation $\R\ni t\to S_t\in Aut(X,\mathscr{B},\mu)$, where $Aut(X,\mathscr{B},\mu)$ denotes the group of measure-preserving automorphisms of $(X,\mathscr{B},\mu)$. Moreover, we assume that this representation is weakly continuous, meaning that for each $f,g\in L^2(X,\mathscr{B},\mu)$ the function $\R\ni t\to \langle f\circ S_t,g\rangle\in \C$ is continuous.\\
\indent Among different classes of flows studied in ergodic theory, one of the most important is the class of horocycle flows acting on the unit tangent spaces of compact surfaces of constant negative curvature. In 1980's, M.\ Ratner \cite{Rat1,Rat3,Rat,Rat2} studied this class  and established several ergodic rigidity phenomena in it. In particular, in \cite{Rat}, she discovered a special way of divergence of orbits (called, in \cite{Rat}, $H_p$-property) for horocycle flows and derived from it important dynamical consequences. The property discovered by Ratner was later named $R_p$-property or Ratner's property (see \cite{Tho}). Ratner's property was further extended by D.\ Witte \cite{Wit} to so called compact Ratner's property and used in the study of the conjugacy problem for unipotent flows.\\
\indent A natural question arises whether there are flows satisfying Ratner's property beyond the class of flows of algebraic origin. This question has been answered positively. Namely, K.\ Frączek and M.\ Lemańczyk \cite{Fr-Lem} have proved that Ratner's property holds for each special flow given by an irrational rotation by $\alpha$, with $\alpha$ having bounded partial qoutients, and the roof function piecewise absolutely continuous\footnote{$f:[0,1)\to\R$ is piecewise absolutely continuous if $[0,1)=[a_0,a_1)\cup[a_1,a_2)\cup...\cup[a_{N-1},a_N)$, $a_0=0,a_N=1$ and $f|_{[a_i,a_{i+1})}$ is absolutely continuous $i=0,...,N-1$.}, satisfying von Neumann's condition (from \cite{Neu}) \begin{equation}\label{1}\int_\T f'd\lambda\neq 0,\end{equation} $\lambda$ stands for Lebesgue measure on $\T$. As a matter of fact, Ratner's property (from \cite{Rat}) was slightly modified in \cite{Fr-Lem} to so called finite Ratner's property, however the dynamical consequences (the rigidity property of joinings, see Theorem \ref{join} below) of the extended property remained valid. We emphasize that when dealing with the aforementioned class of special flows over irrational rotations, in fact, we deal with flows different from (not isomorphic to) horocycle flows. Indeed, horocycle flows are mixing, in fact, they have Lebesgue spectrum of infinite multiplicity (e.g. \cite{Kat-Tho}). It now follows by the classical result by A. Kochergin \cite{Koc} that if we consider special flows over an irrational rotation with the roof function of bounded variation (and such are roof functions from the aforementioned class), then the resulting flows are not mixing. (We note in passing that if we want to be sure that we do not consider, up to isomorphism, horocycle flows, we need to put some smoothness condition on the roof function; otherwise, if we fix an irrational rotation by $\alpha$ and let the roof function be arbitrary continuous then, up to isomorphism, we will obtain every horocycle flow; indeed, every horocycle flow is loosely Bernoulli \cite{Rat3}). More than that, it follows from \cite{Fr-Le} that every special flow over an irrational rotation with the roof function of bounded variation has singular spectrum and is spectrally disjoint from all mixing flows. This shows that flows considered in \cite{Fr-Lem} are, from the dynamical point of view, completely different from horocycle flows.\\
\indent Note that if $f:\T\to\R$ is a function of bounded variation, it has only countably many discontinuity points, say, $0\leq y_1,y_2,...$, all of them are discontinuities of the first order; let $s_i$ denote the jump at $y_i$. In fact, every countable subset of $\T$, in particular dense subsets, can be the set of discontinuities of some $f$ as above. In order to study Ratner's property, we make two introductory remarks. First of all, it seems that finite Ratner's property is too strong to hold in such a general class of flows. That is why we need to weaken this property to so called weak Ratner's property. We borrow the latter notion from \cite{Fra-Lem}, where it was used to study special flows over two-dimensional rotations.\footnote{Dynamical consequences (rigidity of joinings) hold for flows with weak Ratner's property.}  Secondly, note that $\int_\T f'd\lambda$ for piecewise absolutely continuous functions is equal to $S(f)$ which is the (finite) sum of jumps. The number $S(f):=\sum_{i=1}^{+\infty}s_i$ is well defined for each $f$ of bounded variation (the series is absolutely convergent because $f$ is of bounded variation). Therefore, it might seem that $S(f)\neq 0$ (corresponding to (\ref{1}) for piecewise absolutely continuous functions) is a ``working'' condition in our general setup. However, it turns out that the presence of singular component in the Lebesgue decomposition \begin{equation}\label{2}f=f_j+f_a+f_s +S\{\cdot\},\end{equation}
 where $f_j$ is the jump function, $f_a$ is absolutely continuous on $\T$, $f_s$ is singular and continuous on $[0,1]$ (see Section \ref{sec} for details) is another obstacle when studying ergodic properties of such flows. Indeed, it follows from a result from \cite{Iw-Lem-Ma} that for each $\alpha$ with unbounded partial qoutients there exists $f=f_s$, $S(f_s)\neq 0$, and the corresponding special flow is not weakly mixing. More than that, this was strenghtened in \cite{Vol} to have, for each irrational rotation, examples of $f=f_s$, $S(f_s)\neq 0$, isomorphic to the suspension flow.\\
\indent Having all this in mind, we restrict ourselves to study special flows with roof functions for which $f_s=0$, i.e. $f=f_j+f_a+S\{\cdot\}~\footnote{$\{x\}$ denotes the fractional part of $x\in\R$.}$ (note that $f_s=0$ for $f$ piecewise absolutely continuous). Such functions form a Banach space $\mathcal{V}$ with the norm given by $\|f\|= Var(f)+\|f\|_{L^1}$. Now, $S(f)=\int_\T f'd\lambda$ is a linear, continuous functional on $\mathcal{V}$ and the von Neumann condition (\ref{1}) is given by $S(f)=S\neq 0$, where $S$ comes from the decomposition (\ref{2}). It follows that the set $\mathcal{U}$ of functions in $\mathcal{V}$ satisfying (\ref{1}) is open in $\mathcal{V}$, it is the complement of a hyperplane.\\
\indent In order to describe the main results of the paper, recall that von Neumman in \cite{Neu} proved that for an arbitrary irrational $\alpha$ and a piecewise absolutely continuous $f$,~(\ref{1}) implies weak mixing of the corresponding special flow. We will generalize this result by showing (Proposition \ref{weak mixing}) that the special flow over the rotation by an arbitrary irrational $\alpha$ with the roof function in $\mathcal{U}$ is weakly mixing.\\
\indent To investigate stronger properties than weak mixing, we will constantly assume that $\alpha$ has bounded partial quotients. Moreover, there will be some assumptions concerning the rate of convergence of the series of jumps of $f$. This may look quite strange because the set of discontiunuity points can be any countable subset of $\T$. It will be shown however that we can introduce some well-ordering on the set of discontinuity points inherited from a natural well-ordering of the set of absolute values of the jumps. This ordering will be optimal when considering the speed of convergence of the series (see Remark \ref{speed}). One of the main results of the paper is the following (cf. Theorem \ref{ratner}.)
\begin{theorem}\label{ratner}
Let $T:\mathbb{T}\to\mathbb{T}$ be the rotation by an irrational $\alpha$ with bounded partial quotients, i.e. $C:=\sup_i\{a_i\}+1<+\infty$, where $[0;a_1,a_2,...]$ stands for the continued fraction expansion of $\alpha$. Let  $f:\mathbb{T}\to \R_+$ be a bounded away from zero function in $\mathcal{U}$ with the set of jumps $\{d_i\}_{i=1}^{+\infty}$ satisfying for some $j\in \N$
$$\sum_{i>j}|d_i|\leq\frac{|S|}{(2+\theta)(2C+1)((2C+1)^j+1)},$$ for some $\theta>0$. Then the special flow $\mathcal{T}^f$ has the $R(\gamma,P)$ property~\footnote{This property is defined in Section 5.} for every $\gamma>0$, where $P:=(\sgn S)p-D$,  $D:=\{\sum_{i=1}^{+\infty}n_id_i \;|\;0\leq n_i<2C+1 \}$ and $p\in \R$ is such that $(p-\eta,p+\eta)\subset(0,|S|)\setminus(D\cup-D)$ for some $\eta>0$.
\end{theorem}
Theorem 1.1. gives a stability result of weak Ratner's property in the class of piecewise absolutely continuous functions satisfying (\ref{1}) in the space $\mathcal{V}$. Namely, it will follow from Theorem~\ref{ratner} that there is an open set $\mathcal{G}\subset \mathcal{V}$ of functions satisfying weak Ratner's property such that $\mathcal{G}+\R$ is dense in $\mathcal{V}$. In particular, the set of functions satisfying weak Ratner's property has non-empty interior in $\mathcal{V}$ and contains every piecewise absolutely continuous function satisfying~(\ref{1}), see Remark \ref{stab}.\\
\indent The property of mild mixing of a (finite) measure preserving transformation has been introduced by H.\ Furstenberg and B.\ Weiss in \cite{Fu-We}. Recall that a finite measure-preserving transformation is mildly mixing if its Cartesian product with an arbitrary ergodic (finite or infinite not of type I) transformation remains ergodic. It is also proved in \cite{Fu-We} that a finite measure-preserving transformation $T:(X,\mathscr{B},\mu)\to(X,\mathscr{B},\mu)$ is mildly mixing if and only if it has no non-trivial rigid factors, i.e.\ $\liminf_{n\to+\infty}\mu(T^{-n}B\triangle B)>0$ for every $B\in \mathscr{B}$ with $\mu(B)\notin \{0,\mu(X)\}$. The notion of mild mixing has been extensively studied by many authors, see e.g.\ \cite{Aa-Li-We,Fu,Lem-Les,Lem-Pa,Sch}. Similar definition can be introduced and similar results hold in case of flows. It is immediate from the definition that the strong mixing property of a flow implies its mild mixing which implies weak mixing property. It follows from \cite{Fr-Lem} that for flows having weak Ratner's property, to prove mild mixing, it suffices to show the absence of partial rigidity (see Sections 2 and 3 for needed definitions).\\
\indent The second main result of the paper is Theorem \ref{rig}, in which the absence of partial rigidity for special flows under the functions from $\mathcal{U}$ and with some assumptions concerning the rate of convergence of the series of jumps is proved. Contrary to \cite{Fr-Lem}, where absence of partial rigidity is proved for all special flows under piecewise absolutely continuous functions and arbitrary irrational $\alpha$, Theorem \ref{rig} is shown to hold for $\alpha$ having bounded partial qoutients.\\
\indent The author would like to thank Professor Mariusz Lemańczyk for many fruitful discussions and suggestions concerning this paper.

\section{Basic notions}

Let $\mathcal{S}=(S_t)_{t\in\R}$ be an ergodic (measurable) flow  on a standard probability space $(X,\mathscr{B},\mu)$. We recall that $\mathcal{S}$ is called {\em mixing} if
$$\lim_{t\to+\infty}\mu(S_tA\cap B)=\mu(A)\mu(B)\;\hbox{for all}\; A,B\in\mathscr{B}.$$
If 
$$\lim_{T\to+\infty}\frac{1}{T}\int_0^T|\mu(S_tA\cap B)-\mu(A)\mu(B)|dt=0,$$
for all $A,B\in\mathscr{B}$ then $(S_t)_{t\in\R}$ is called {\em weakly mixing}.
Recall also that the weak mixing of $\mathcal{S}$ is equivalent to the ergodicity of $(S_t\times S'_t)_{t\in\R}$ acting on $(X\times X',\mathscr{B}\otimes\mathscr{B'},\mu\times\mu')$ for each ergodic $(S'_t)_{t\in\R}$ acting on a standard probability space $(X',\mathscr{B'},\mu')$. It easily follows that mixing implies mild mixing which in turn implies weak mixing.\\
\indent Assume that $\mathscr{A}\subset\mathscr{B}$ is a {\em factor} of $(S_t)_{t\in\R}$, i.e.\ $\mathscr{A}$ is an $(S_t)_{t\in\R}$-invariant sub-$\sigma$-algebra. Let $(t_n)_{n\in\N}$ be a sequence of real numbers such that $t_n\to+\infty$. We say that the factor $\mathscr{A}$ is {\em rigid} along $(t_n)$ if
$$\lim_{n\to+\infty}\mu(A\cap S_{-t_n}A)= \mu(A)$$
for every $A\in\mathscr{A}$. In particular, $(S_t)$ is rigid along $(t_n)$ if $\mathscr{A}=\mathscr{B}$. A flow $(S_t)_{t\in\R}$ is called {\em partially rigid} along $(t_n)$ if there exists $0<u\leq 1$ such that
$$\liminf_{n\to+\infty} \mu(B \cap S_{-t_n}B)\geq u\mu(B)\;\; \hbox{for every}\; B\in\mathscr{B}.$$
\indent Assume that $T$ is an ergodic automorphism on $(X,\mathscr{B},\mu)$. A measurable function $f:X\to\R$ determines a cocycle $f^{(\cdot)}(\cdot):\Z\times X\to \R$ given by
$$f^{(m)}(x)=\left\{\begin{array}{cll} f(x)+f(Tx)+...+f(T^{m-1}x)& \text{if}& m>0 \\
0&\text{if}& m=0\\
 -(f(T^mx)+...+f(T^{-1}x))&\text{if}&m<0.\end{array}\right.$$

Denote by $\lambda$ Lebesgue measure on $\R$. If $f:X\to\R$ is a strictly positive $L^1$ function, then by $\mathcal{T}^f=(T_t^f)_{t\in\R}$ we will mean the corresponding special flow under $f$ acting on $(X^f,\mathscr{B}^f,\mu^f)$, where $X^f:=\{(x,s)\in X\times\R:\; 0\leq s<f(x)\}$, and $\mathcal{B}^f$ ($\mu^f$) is the restriction of $\mathscr{B}\otimes \mathscr{B}(\R)$ ($\mu\times \lambda$) to $X^f$. Under the action of the flow $\mathcal{T}^f$ each point in $X^f$ moves vertically with unit speed and we identify the point $(x,f(x))$ with $(Tx,0)$. More precisely, if $(x,s)\in X^f$ then
$$T_t^f(x,s)=(T^nx,s+t-f^{(n)}(x)),$$
where $n\in\Z$ is unique such that
$$f^{(n)}(x)\leq s+t<f^{(n+1)}(x).$$ 
 It is well-known (due to \cite{Neu}) that the special flow $(T_t^f)$ is weakly mixing if and only if for every $s\in\R\setminus\{0\}$ the equation
$$\frac{\psi(Tx)}{\psi(x)}=e^{2\pi isf(x)}$$
has no measurable solution $\psi:X\to \mathbb{S}^1:=\{z\in\C:\: |z|=1\}$.\\

\section{Joinings}
Let $\mathcal{S}=(S_t)_{t\in\R}$, $\mathcal{T}=(T_t)_{t\in\R}$ be two ergodic flows defined on $(X,\mathscr{B},\mu)$ and $(Y,\mathscr{C},\nu)$ respectively. By a {\em joining} between $\mathcal{S}$ and $\mathcal{T}$ we mean any probability $(S_t\times T_t)_{t\in\R}$-invariant measure on $(X\times Y,\mathscr{B}\otimes \mathscr{C})$ whose projections on $X$ and $Y$ are equal to $\mu$ and $\nu$ respectively. The set of joinings between $\mathcal{S}$ and $\mathcal{T}$ is denoted by $J(\mathcal{S},\mathcal{T})$. The subset of ergodic joinings is denoted by $J^e(\mathcal{S},\mathcal{T})$. Ergodic joinings are exactly extremal points in the simplex $J(\mathcal{S},\mathcal{T})$. Let $\{A_n:\:n\in\N\}$  and $\{B_n:\:n\in\N\}$ be two countable families in $\mathscr{B}$ and $\mathscr{C}$ respectively, which are dense for the pseudo-metrics $d_\mu(A,A')=\mu(A\triangle A')$ and $d_\nu(B,B')=\nu(B\triangle B')$ respectively. Let us consider the metric $d$ on $J(\mathcal{S},\mathcal{T})$ defined by
$$d(\rho,\rho')=\sum_{m,n\in\N}\frac{1}{2^{m+n}}|\rho(A_n\times B_m)-\rho'(A_n\times B_m)|,\;\rho,\rho'\in J(\mathcal{S},\mathcal{T}).$$
Endowed with the topology corresponding to $d$, which we will refer to as the weak topology, the set $J(\mathcal{S},\mathcal{T})$ is compact.\\

 Denote by $\mu\times_\mathscr{A}\mu\in J(\mathcal{S},\mathcal{S})$ the relatively independent joining of the measure $\mu$ over the factor $\mathscr{A}$, i.e. 
$$(\mu\times_\mathscr{A}\mu)(D)=\int_{X/\mathscr{A}}(\mu_{\bar{x}}\times\mu_{\bar{x}})(D)\,d\bar{\mu}(\bar{x})$$
for $D\in \mathscr{B}\otimes \mathscr{B}$, where $\{\mu_{\bar{x}}:\;\bar{x}\in X/\mathscr{A}\}$ is the disintegration of the measure $\mu$ over the image $\bar{\mu}$ of $\mu$ via the factor map $\Pi:X\to X/\mathscr{A}$.\\

\indent For every $t\in\R$, by $\mu_{S_t}\in J^e(\mathcal{S},\mathcal{S})$, we will denote the graph joining determined by $(S_t)_{t\in\R}$, i.e. $\mu_{S_t}(A\times B)=\mu(A\cap S_{-t}B)$ for $A,B\in\mathscr{B}$. Then $\mu_{S_t}$ is concentrated on the graph of $S_t$.\\

Recall that in general the notions of (absence of) partial rigidity and mild mixing are not related. We have however the following.
\begin{lemma} {\em \cite{Fr-Lem}}\label{mild mixing}
Let $\mathcal{S}$ be an ergodic flow on $(X,\mathscr{B},\mu)$ which is a finite extension of each of its non-trivial factors. If the flow  $\mathcal{S}$ is not partially rigid then it is mildly mixing. 
 \end{lemma}
\section{Continued fraction expansion of an irrational number $\alpha$.}
We denote by $\T$ the circle group $\R/\Z$ which will be identified with the interval $[0,1)$ with addition $mod\, 1$. For a real number $t$ denote by $\{t\}$ its fractional part and by $\|t\|$ its distance to the nearest integer number. For an irrational $\alpha\in\T$ denote by $(q_n)_{n=0}^{+\infty}$ its sequence of denominators, that is, we have 
$$\frac{1}{2q_nq_{n+1}}<\left|\alpha-\frac{p_n}{q_n}\right|<\frac{1}{q_nq_{n+1}},$$
where $q_0=1,\; q_1=a_1,\;\; q_{n+1}=a_{n+1}q_n+q_{n-1}$, $p_0=0,\; p_1=1,\;\; p_{n+1}=a_{n+1}p_n+p_{n-1}$ and $[0;a_1,a_2,...]$ stands for the continued fraction expansion of $\alpha$. One says that $\alpha$ has {\em bounded partial quotients} if the sequence $(a_n)_{n=1}^{+\infty}$ is bounded. In this case, if we set $C:=\sup\{a_n:\; n\in\N\}+1$ then $q_{n+1}\leq Cq_n$ and
$$\frac{1}{2Cq_n}\leq\frac{1}{2q_{n+1}}<\|q_n\alpha\|<\frac{1}{q_{n+1}}<\frac{1}{q_n}$$
for each $n\in\N$. The following lemma is well-known.
\begin{lemma}\label{dist} Let $\alpha\in\T$ be irrational with bounded partial quotients. Then there exist positive constants $C_1,C_2$ such that for every $k\in\N$ the lenghts of intervals $J_1,...,J_k$ arisen from the partition of $\T$ by $0,-\alpha,...,-(k-1)\alpha$ satisfy $\frac{C_2}{k}\leq|J_j|<\frac{C_1}{k}$ for each $j=1,...,k$. 
\end{lemma}

\label{weak}\section{Weak Ratner's property}
In this section we recall  the notion of weak Ratner's property introduced in \cite{Fra-Lem} and we list results from \cite{Fra-Lem} needed in what follows.
\begin{definition}{\em \cite{Fra-Lem} Let $(X,d)$ be a $\sigma$-compact metric space, $\mathscr{B}$ the $\sigma$-algebra of Borel subsets of $X$, $\mu$ a Borel probability measure 
on $(X,d)$ and let $\mathcal{S}=(S_t)_{t\in\R}$ be a flow on $(X,\mathscr{B},\mu)$. Let $P\subset \R\setminus\{0\}$ be a compact subset and $t_0\in\R\setminus\{0\}$. The flow $(S_t)_{t\in\R}$ is said to have the property $R(t_0,P)$ if for every $\epsilon>0$ and $n\in\N$ there exist $\kappa=\kappa(\epsilon)$, $\delta=\delta(\epsilon,N)>0$ and a subset $Z=Z(\epsilon,N)\in\mathscr{B}$ with $\mu(Z)>1-\epsilon$ such that if $x,x'\in Z$, $x'$ is not in the orbit of $x$, and $d(x,x')<\delta$, then there are $M=M(x,x')\geq N$, $L=L(x,x')\geq N$ such that $\frac{L}{M}\geq \kappa$ and there exists $\rho=\rho(x,x')\in P$ such that 
$$\frac{1}{L}\big|\{n\in\Z\cap[M,M+L] :\; d(S_{nt_0}(x),S_{nt_0+\rho}(x'))<\epsilon\}\big|>1-\epsilon.$$}
\end{definition}
Moreover, we say that $(S_t)_{t\in\R}$ has  {\em weak Ratner's property}  or {\em $R(P)$-property} if the set of $s\in\R$ such that the flow $(S_t)_{t\in\R}$ has the $R(s,P)$-property is uncountable.\\
\indent We will constantly assume that $(S_t)_{t\in\R}$ satisfies the following ``almost continuity'' condition: for every $\epsilon>0$ there exists $X(\epsilon)\in\mathscr{B}$ with $\mu(X(\epsilon))>1-\epsilon$ such that for every $\epsilon'>0$ there exists $\epsilon_1>0$ such that 
$d(S_tx,S_{t'}x)<\epsilon'$ for all $x\in X(\epsilon)$ and $t,t'\in[-\epsilon_1,\epsilon_1]$.\\
\indent Note that if $(T_t^f)_{t\in\R}$ is a special flow acting on $(X^f,\mathscr{B}^f,\mu^f)$ equipped with the metric $d_1((x,t),(y,s))=d(x,y)+|t-s|$, then the above condition holds.

\begin{theorem}{\em \cite{Fra-Lem}} \label{join}
Let (X,d) be a $\sigma$-compact metric space, $\mathscr{B}$ the $\sigma$-algebra of Borel subsets of X and $\mu$ a probability Borel measure on (X,d). Let $(S_t)_{t\in\R}$ be a weakly mixing flow on the space $(X,\mathscr{B},\mu)$ that satisfies the $R(P)$-property, where $P\subset\R\setminus\{0\}$ is a nonempty compact set. Assume that $(S_t)_{t\in\R}$ satisfies the ``almost continuity'' condition. Let $(T_t)_{t\in\R}$ be an ergodic flow on $(Y,\mathscr{C},\nu)$  and let $\rho$ be an ergodic joining of $(S_t)_{t\in\R}$ $(T_t)_{t\in\R}$. Then either $\rho=\mu\times \nu$ or $\rho$ is a finite extension of $\nu$.
\end{theorem}
\begin{remark} { \em Let $\mathcal{S}$ be an ergodic flow on $(X,\mathscr{B},\mu)$. Assume that for each ergodic flow $\mathcal{T}$ acting on  $(Y,\mathscr{C},\nu)$ an arbitrary ergodic joining $\rho$ is either the product measure or $\rho$ is a finite extension of $\nu$. Then $\mathcal{S}$ is a finite extension of each of its non-trivial factors.  Indeed, suppose that $\mathscr{A}\subset \mathscr{B}$ is a non-trivial factor. Let us consider the action of $\mathcal{S}|_{\mathscr{A}}$ on $(X/\mathscr{A},\mathscr{A},\overline{\mu})$ and the natural joining $\mu|_{\mathscr{A}}\in J^e(\mathcal{S},\mathcal{S}|_{\mathscr{A}})$ determined by $\mu|_{\mathscr{A}}(B\times A)=\mu(B\cap A)$ for all $B\in \mathscr{B}$, $A\in \mathscr{A}$. Clearly the action $\mathcal{S}\times (\mathcal{S}|_{\mathscr{A}})$ on $(X\times (X/\mathscr{A}),\mathscr{B}\otimes \mathscr{A},\mu|_{\mathscr{A}})$ is isomorphic to the action of $\mathcal{S}$ (via the projection on $(X,\mathscr{B},\mu)$). Since the measure $\mu|_{\mathscr{A}}$ is not the product measure, the action  $\mathcal{S}\times (\mathcal{S}|_{\mathscr{A}})$ is a finite extension of $\mathcal{S}|_\mathscr{A}$.
}
\end{remark}
\begin{proposition}\label{condi}{\em \cite{Fra-Lem}} Let (X,d) be a $\sigma$-compact metric space, $\mathscr{B}$ be the $\sigma$-algebra of Borel subsets of $X$, $\mu$ a probability Borel measure on $(X,d)$. Assume that $T:(X,\mathscr{B},\mu)\to(X,\mathscr{B},\mu)$ is an ergodic isometry and $f:X\to\R$ is a bounded positive measurable function which is bounded away from zero. Let $P\subset\R\setminus\{0\}$ be a compact set. Assume that for every $\epsilon>0$ and $N\in\N$ there exist $\kappa=\kappa(\epsilon)>0$, $\delta=\delta(\epsilon,N)>0$ and a subset $Z=Z(\epsilon,N)\in\mathscr{B}$ with $\mu(Z)>1-\epsilon$ such that if $x,x'\in Z$, and $d(x,x')<\delta$, then there are $M=M(x,x')\geq N$, $L=L(x,x')\geq N$ such that $\frac{L}{M}\geq \kappa$ and there exists $p=p(x,x')\in P$ such that 
$$\frac{1}{L}\left|\{n\in\Z\cap[M,M+L] :\; |f^{(n)}(x)-f^{(n)}(x')-p|<\epsilon\}\right|>1-\epsilon.$$
Suppose that $\gamma\in \R$ is a positive number such that the $\gamma$-time automorphism $T^f_\gamma:X^f\to X^f$ is ergodic. Then the special flow $\mathcal{T}^f$ has the $R(\gamma,P)$-property.
\end{proposition}
\subsection{Properties of bounded variation functions without singular continuous component in the Lebesgue decomposition} \label{sec}

Let $T:\T\to\T$ be an irrational rotation by $\alpha$ with the sequence of denominators $(q_n)_{n=1}^{+\infty}$ and $f:\mathbb{T} \to \R$ a function with bounded variation $(f\in BV(\T))$\footnote{We constantly assume that all considered functions of bounded variation are left-continuous.}. It is well known that $f\in L^1(\T)$ and the following inequality holds.
\begin{proposition}\label{koks}{\em(Denjoy-Koksma inequality; see e.g.\cite{Co-Fo-Si,Kui-Nid})} If $f:\T\to\R$ is a function of bounded variation then
$$\left|\sum_{k=0}^{q_n-1}f(x+k\alpha)-q_n\int_\T f\,d\lambda\right|\leq Varf,$$
for every $x\in \T$ and $n\in \N$.
\end{proposition}
 Let us identify the function $f:\T\to \R$, with the function $\bar{f}:[0,1]\to \R$, by $\bar{f}(x)=f(x)$ if $x\in[0,1)$, and $\bar{f}(1)=\lim_{a\to 1^-}f(a)$ which exists since $f\in BV(\T)$. Then $\bar{f}\in BV([0,1])$. (Note that if  $f$ is continuous from the left, so is $\bar{f}$.) By the Lebesgue decomposition, see e.g. \cite{Fre}, it follows that
$$\bar{f}=\bar{f_a}+\bar{f_j}+\bar{f_s},$$
where $\bar{f_j}(x)=\sum_{y<x}(\bar{f}(y^+)-\bar{f}(y^-))$~\footnote{$f(y^{\pm}):=\lim_{z\to y^{\pm}}f(z)$.} is the jump function. Note that $\bar{f_j}$ has countably many discontinuity points, say $\beta_2,\beta_3,\ldots$, and the series of jumps is absolutely convergent (since $\bar{f}\in BV([0,1])$). Furthermore, $\bar{f_s}$ is singular and continuous, and $\bar{f_a}$ absolutely continuous on $[0,1]$. This decomposition is unique up to constants. It follows that we can decompose $f$ as
\begin{equation}\label{eq} f=f_a+f_j+f_s+S\{\cdot\},\end{equation}
where $f_j:=\bar{f}_j|_{[0,1)}$ is the jump function on $\T$ with the set of jumps $\{d_i\}_{i=1}^{+\infty}$ at $\{\beta_i\}_{i=1}^{+\infty}$ respectively ($\beta_1=0$ and $d_1$ can be equal to zero), $f_s:=\bar{f}_s|_{[0,1)}$ is singular, continuous on $\T\setminus\{0\}$ (the only discontinuity point can be 0), and $f_a:=\bar{f}_a|_{[0,1)}-S\{\cdot\}$ is absolutely continuous on $\T$ ($S=\bar{f}_a(1)-\bar{f}_a(0)$). Recall that $BV(\T)$ is a Banach space with the norm $\|f\|:=Varf+ \|f\|_{L^1}$. Let us set $BV_{a+j}(\T):=\{f\in BV(\T):\: f_s=0\}(=\mathcal{V})$ which is a Banach space (it is a closed subspace of $BV(\T)$). For the functions in $BV_{a+j}(\T)$ we have that $S:=\sum_{i=1}^{+\infty}d_i$ is just the sum of the jumps of the function $f$ and $Varf=\sum_{i=1}^{+\infty}|d_i|+ \|f_a'\|_{L^1}+2S$. Consider the set $\mathcal{U}:=\{f\in BV_{a+j}(\T):\:S\neq 0 \}$. Notice that all piecewise absolutely continuous functions satisfying (\ref{1}), which from now on will be called von Neumann's functions, belong to $\mathcal{U}$. Von Neumann's functions are precisely all functions in $BV_{a+j}(\T)$ for which there are only finitely many discontinuity points. Moreover, $\mathcal{U}$ is an open set in $BV_{a+j}(\T)$. In fact, it is the complement of a hyperplane.
\begin{proposition}\label{dense} The set of von Neumann's functions is dense in $\mathcal{U}$.
\end{proposition}
Proof: Let $f\in \mathcal{U}$ be a function with the set of  discontinuities $\{\beta_i\}_{i=1}^{+\infty}$ and the set of corresponding jumps $\{d_i\}_{i=1}^{+\infty}$. We will construct a sequence $(f_n)_{n=1}^{+\infty}\subset \mathcal{U}$ of von Neumann's functions which tends to $f$ as $n\to +\infty$. The series $\sum_{i=1}^{+\infty}|d_i|$ is convergent, so there exists $j_n\in \N$ such that $\sum_{i>j_n}|d_i|<\frac{1}{3n}$. Let us consider the discontinuity points $\beta_i$, ${i=1,...,j_n}$. By permuting them, write $0\leq\beta_{k_1}<\beta_{k_2}<...<\beta_{k_{j_n}}$ where $\{1,...,j_n\}=\{k_1,...,k_{j_n}\}$. Let $g_{j_n}(x)=\sum_{i\leq k} d_i$ if $\beta_{k_r}\leq x<\beta_{k_{r+1}}$, $r=1,...,j_n$, $\beta_{j_{n+1}}=\beta_{k_1}$ and set $S_n:=\sum_{i=1}^{j_n}d_i$. Then let 
$$f_n(x):=f_a(x)+g_{j_n}(x)+S_n\{x\}\in BV_{a+j}(\T).$$
It follows that the only points of discontinuity of $f_n$ are $\beta_{k_1},...,\beta_{k_{j_n}}$ (perhaps $f_n$ is continuous at $\beta_{k_1}$). Moreover, $|S-S_n|<\frac{1}{3n}$. Since $S\neq 0$, there exists $n_0\in\N$ such that $S_n\neq 0$ if $n\geq n_0$. Thus, $f_n$ is a von Neumann's function. Now, consider $(f-f_n)(x)=(f_j-g_{j_n})(x)+(S-S_n)\{x\}$. We have
$$Var(f-f_n)\leq Var((S-S_n)\{\cdot\})+Var((f_j-g_{j_n}))=2|S-S_n|+\sum_{i>j_n}|d_i|\leq \frac{1}{n},$$
so $f_n\to f$ in $BV_{a+j}(\T)$ (obviously $f_n\to f$ in $L^1$).  \hfill  $\square$ 

\vspace{2ex}

As an arbitrary $f\in\mathcal{U}$ is left-continuous, 
\begin{equation}\label{fac}
f:=f_{ac}+f_{pl},
\end{equation} 
 where $f_{ac}:\mathbb{T}\to \R$ is absolutely continuous with zero mean and $f_{pl}:\mathbb{T}\to \R$ belongs to $\mathcal{U}$ and  $f_{pl}'(x)=S$ for all $x\in\T\setminus\{\beta_i\}_{i=1}^{+\infty}$. Indeed, we set $$f_{pl}(x):=\sum_{i=1}^{+\infty}d_i\{x-\beta_i\}+c,$$
for some constant $c\in \R$ ($c:=\int_\T f_a(x)\:d\lambda$, cf.\ (\ref{2})). Note that the points of discontinuity and the size of jumps of $f$ and $f_{pl}$ are the same. Indeed, it is an easy consequence of the decomposition (\ref{2}) ($f_{pl}(x)=f_j(x)+S\{x\}+c$). \\
\indent We consider $f\in BV_{a+j}(\T)$. Then $\int_\T f'd\lambda=S$. Moreover, for each interval $[a,b]\subset [0,1)$ we have
\begin{equation} \label{inpa}
\int_a^bf'(x)d\lambda=f(b^-)-f(a^+)-\sum_{\{i;\;\beta_i\in [a,b]\}}d_i=f(b^-)-f(a^+)-\sum_{x\in[a,b]}f_j(x^+)-f_j(x).
\end{equation}
Indeed, 
 $$\int_a^bf'd\lambda=\int_a^bf_a'd\lambda+\int_a^bS\;d\lambda=f_a(b^-)-f_a(a^+)+S(b-a)=$$
$$f_a(b^-)-f_a(a^+)+S(b-a)+f_j(b^-)-f_j(a^+) -(f_j(b^-)-f_j(a^+))=$$
$$f_a(b^-)-f_a(a^+)-\sum_{\{i;\;\beta_i\in [a,b]\}}d_i.$$

\subsection{Von Neumann's functions considered up to cohomology do not exhaust $\mathcal{U}$}
\indent We say that $f:\T\to \R$ is a {\em coboundary} if there exists a measurable function  $j:\T\to\R$ such that $f(x)=j(x)-j(Tx)$ for almost every $x\in\T$. One says that $f$ and $g$ are {\em cohomologous} if their difference is a coboundary. In the previous section, we have shown that von Neumann's functions are dense in $\mathcal{U}$, whence a natural question arises whether there are some functions $f\in \mathcal{U}$ which are not cohomologous to any von Neumann's function.\\
\indent To answer positively this question, let us consider a function $f\in \mathcal{U}$ with an infinite set of discontinuities $\{\beta_i\}_{i=1}^{+\infty}\subset \Q$ such that $\beta_i-\beta_j\notin \Z+\Z\alpha$ whenever $i\neq j$. Let the jump $d_i$ at $\beta_i$ be positive for each $i\in \N$. Moreover, we assume that for every $\epsilon>0$ there exists $N_\epsilon$ such that \begin{equation}\label{coh}\sum_{i>N_\epsilon}d_i\leq \epsilon\min\{d_1,...,d_{N_\epsilon}\}.\end{equation} 
Take $g$ a von Neumann's function with the set of  discontinuities  $\{\gamma_j\}_{j=1}^k$ and the corresponding set of jumps $\{m_j\}_{j=1}^k$, $S(g)=\int_\T g'd\lambda=\sum_{j=1}^km_j\neq 0$, and suppose that $(f-g)(x)=j(x)-j(Tx)$ for some measurable function $j:\T\to \R$. Let $(q_n)$ be the sequence of denominators of $\alpha$. Then $(f-g)^{(q_n)}\to 0$ in measure as $n\to +\infty$. We must have that $S(f)=S(g)$, if not $f-g\in \mathcal{U}$ and by Theorem \ref{weak mixing} below, we get a contradiction. By (\ref{2}), $f-g=(f_j-g_j)+(f_a-g_a)=(f_j-g_j+c)+(f_a-g_a-c)$\footnote{$c:=\int_\T (f_a-g_a)d\lambda$.}. Using Denjoy-Koksma inequality (Proposition \ref{koks}), $(f_a-g_a-c)^{(q_n)}\to 0$
uniformly as $n\to +\infty$. It follows that $(f_j-g_j+c)^{(q_n)}\to 0$ in measure as $n\to +\infty$, and we will show that this is impossible under our additional assumption (\ref{coh}). 
\begin{remark} \label{Co-Pi} {(e.g.\ \em\cite{Co-Pie})  Without loss of generality we can assume that the difference between any two discontinuity points of $f_j-g_j$ is not a multiple of $\alpha$. Indeed, if $\delta-\delta'\in \Z+\Z\alpha$, $\chi_{[\delta,\delta')}(x)-(\delta'-\delta)$ is a coboundary.}
\end{remark}
Since the number of discontinuities of $f$ is infinite and $\beta_i-\beta_j\notin \Z+\Z\alpha$, after applying Remark \ref{Co-Pi}, there are still infinitely many points of discontinuity of $f_j-g_j$ left and the difference between any two of them is not a multiple of $\alpha$. Without loss of generality, we can assume that the set of points of discontinuity is $\{\beta_i\}_{i=1}^{+\infty}\cup\{\gamma_j\}_{j=1}^k$.   

\begin{lemma}{\em(\cite{Fr-Lem-Les}, Lemma 2.3.)} Let $\alpha$ be irrational with bounded partial qoutients and let $\beta\in (\mathbb{Q}+\mathbb{Q}\alpha)\setminus (\Z +\Z\alpha)$. Then there exists $c>0$ such that for each $m\in \N$ the length of each interval in the partition of $\T$ arisen from $0,-\alpha,...,-(m-1)\alpha,\beta,\beta-\alpha,...,\beta-(m-1)\alpha$ is at least $\frac{c}{m}$. 
\end{lemma}
 \indent Let us fix $\epsilon_0>0$ such that $N_{\epsilon_0}\geq 4k$ and, for $m\in\N$, let $\mathcal{I}_m$ be the partition of the circle given by the points $\beta_i-j\alpha$, $i=1,...,N_{\epsilon_0}$, $j=0,...,m-1$. 
\begin{remark} \label{leng} {\em By the above lemma, there exists a constant $c':=c'(\beta_1,...,\beta_{N_{\epsilon_0}})>0$ such that the length of each interval in partition $\mathcal{I}_m$ is at least $\frac{c'}{m}$. It also follows by Lemma \ref{dist} that the length of each interval in $\mathcal{I}_m$ is at most $\frac{C_1}{m}$.~\footnote{If $m=q_n$ for some $n\in \N$, $C_1=1$.}}
\end{remark}
Let us consider the points $\gamma_i-s\alpha$, $i=1,...,k$, $s=0,...,q_n-1$, the number of such points is $kq_n$. The number of intervals in $\mathcal{I}_{q_n}$ is at least $4kq_n$ (because $N_{\epsilon_0}\geq 4k$).  Fix $1\leq i\leq N_{\epsilon_0}$. For $s=0,...,q_n-1$ let $I_{i,s}, I'_{i,s}\in \mathcal{I}_{q_n}$ be the intervals with the right endpoint $\beta_i-s\alpha$ and the left endpoint $\beta_i-s\alpha$ respectively. Then there exists $1\leq i'\leq N_{\epsilon_0}$ such that $|A_{q_n,i'}|>q_n/2$, where
$$A_{q_n,i'}:=\left\{0\leq s\leq q_n-1:\;\left(\bigcup_{j=1}^k\{\gamma_j-r\alpha\}_{r=0,...,q_n-1}\right)\cap\left(I_{i',s}\cup I'_{i',s}\right)=\emptyset \right\}.$$

Set $\mathcal{H}_{q_n}(i'):=\bigcup_{s\in A_{q_n,i'}}I'_{i',s}$ and $\mathcal{H'}_{q_n}(i'):=\bigcup_{s\in A_{q_n,i'}}I_{i',s}$. By definition \begin{equation}\label{dis}\left(\bigcup_{s=1}^{q_n-1}\{\gamma_i-s\alpha\}_{i=1}^k\right)\cap\mathcal{H}_{q_n}(i')=\emptyset,\;\;\;\;
\left(\bigcup_{s=1}^{q_n-1}\{\gamma_i-s\alpha\}_{i=1}^k\right)\cap\mathcal{H'}_{q_n}(i')=\emptyset.\end{equation}

It follows that $|\mathcal{H}_{q_n}(i')|, |\mathcal{H'}_{q_n}(i')|\geq\frac{q_n}{2}$ (each interval in the above families has length at least $\frac{c'}{q_n}$). Let us fix $I_{i',s}\in\mathcal{H}_{q_n}(i')$ and $I'_{i',s}\in\mathcal{H'}_{q_n}(i')$. Take any $x\in I_{i',s}, y\in I'_{i',s}$. By Remark \ref{leng}, the length of the interval $[x,y]$ is at most $2/q_n$. Moreover, the only points of discontinuity of $(f-g)^{(q_n)}$ in $[x,y)$ are $\beta_{i'}-s\alpha$ and some of the form $\beta_i-r\alpha$, $i>N_\epsilon$,~$r=0,...,q_n-1$. Let us fix $i_1>N_\epsilon$ and consider the points $\beta_{i_1}-t\alpha$,~$\beta_{i_1}-r\alpha$ for some $t\neq r$. Then
$$\|(\beta_{i_1}-t\alpha)-(\beta_{i_1}-r\alpha)\|\geq\|q_{n-1}\alpha\|\geq 1/2q_n.$$
It follows that the number of discontinuities of $(f-g)^{(q_n)}$ in $[x,y)$ of the form $\beta_{i_1}-r\alpha$ for some $r=0,...,q_n-1$ is at most 5. Thus 
\beab(f-g)^{(q_n)}(y)-(f-g)^{(q_n)}(x)=\sum_{i=1}^{+\infty}\left(\sum_{r=0}^{q_n-1}\chi_{[x,y)}(\{\beta_i-r\alpha\})\right)d_i + \\ \sum_{j=1}^{k}\left(\sum_{r=0}^{q_n-1}\chi_{[x,y)}(\{\gamma_j-r\alpha\})\right)m_j= d_{i'}+ \sum_{i>N_{\epsilon_0}}n_id_i,\eeab
for some integers $|n_i|\leq 5$. We get that $(f-g)^{(q_n)}(y)-(f-g)^{(q_n)}(x)\geq (1-5\epsilon_0)d_{i'}$, which is a contradiction with $((f-g)^{(q_n)})_*\lambda\to \delta_{0}$ as $n\to +\infty$.

\section{Weak mixing of special flows when the roof function is in $\mathcal{U}$}
Let $T$ be an ergodic automorphism acting on $(X,\mathscr{B},\mu)$ and $f:X\to \R_+$ in $L^1(X,\mathscr{B},\mu)$. Assume that $T$ is rigid and let $(q_n)_{n=1}^{+\infty}$ be a rigidity sequence. We will state a criterion for weak mixing of special flows over rigid systems.
\begin{lemma}\label{crit}{\em(cf. \cite{Fra-Lem}, Proposition 2.1.)} Under the above assumptions, assume additionally that there exists $0<c<1$ such that
$$\limsup_{n\to+\infty}\left|\int_Xe^{2\pi irf^{(q_n)}(x)}\,d\mu(x)\right|<c,$$
for all $|r|\in \R$ large enough. Then the special flow $(T_t^f)_{t\in\R}$ is weakly mixing.
\end{lemma}
Proof: Suppose on the contrary that for some $s\neq 0$ and a measurable $\psi:X\to \mathbb{S}^1$
$$\frac{\psi(Tx)}{\psi(x)}=e^{2\pi isf(x)}.$$
Then for all $k\in \Z\setminus\{0\}$ we have
$$\left|\int_X(\psi(T^{q_n}x)\overline{\psi(x)})^k\,d\mu(x)\right|=\left|\int_Xe^{2\pi iksf^{(q_n)}(x)}\,d\mu(x)\right|<c<1,$$
for $n,k$ large enough. Since clearly $(\psi\circ T^{q_n})^k\overline{\psi^k}\to 1$ in measure as $n\to +\infty$, we obtain a contradiction. \hfill $\square$  

\vspace{2ex}

We will now prove that all special flows under the roof functions from $\mathcal{U}$ are weakly mixing. The proof of this fact is based on the proof of Theorem 3 from \cite{Lem-Pa} concerning the ergodicity of real cocycles over irrational rotations.
\begin{proposition}\label{weak mixing} Let $T:\T\to\T$ be an arbitrary irrational rotation by $\alpha$.  
Let $f:\T\to\R_+$, $f\in\mathcal{U}$. Then the special flow $\mathcal{T}^f$ is weakly mixing. 
\end{proposition}
Proof: By Theorem \ref{dense}, there exists a sequence $(f_n)_{n=1}^{+\infty}$ of von Neumann's functions which tends to $f$ in $BV_{a+j}(\T)$. It follows that there exists $n_0$ such that for $n\geq n_0$, $|S(f_n)|\neq 0$ and $Var (f-f_n)<|S(f_n)|$ (this holds because $Var(f-f_n)\to 0$ and $S(f_n)\to S(f)\neq 0$). Let $g:=f_{n_0+1}$ and $S:=|S(g)|$. Then $g$ is a von Neumann's function and let $K$ denote its number of discontinuities. We will prove that if $c$ is any constant satisfying $Var(f-g)/S<c<1$, then
\begin{equation}\label{lim}\limsup_{n\to+\infty}\left|\int_\T e^{2\pi irf^{(q_n)}(x)}d\lambda(x)\right|<c\end{equation}
for every $|r|$ large enough. This will show weak mixing in view of Lemma \ref{crit}.\\
\indent Recall (see (\ref{fac})) that $f(x)=f_{pl}(x)+f_{ac}(x)$, where $f_{ac}:\T\to \R$ is absolutely continuous with zero mean and $f_{pl}(x)= f_j(x)+S(f)\{x\}+c'$ ($c'=\int_\T f_{ac}(x)\,d\lambda(x)$). It follows from Denjoy-Koksma inequality that $f_{ac}^{(q_n)}\to 0$ uniformly as $n\to +\infty$, so to prove~(\ref{lim}), without loss of generality, assume that $f=f_{pl}$.  By the proof of Proposition \ref{dense}, $g$ is piecewise linear ($g(x)=g_j(x)+S\{x\}+c'$). Let $h:=f-g$, then, by definition, $Var(h)<S$. Let 
$$I_{r,q}=\int_{\T}e^{2\pi ir(g^{(q)}(x)+h^{(q)}(x))}d\lambda(x)$$
for $r\neq 0$ and $q\geq 1$. \\
\indent Since $g$ is piecewise linear and $g'=S$, we get $\frac{1}{q}|g'^{(q)}|= |S|$ so for all $q$ we have
$$|g'^{(q)}(x)|=|S|q,$$
for almost every $x\in\T$.\\
\indent Denote by $x_1\leq...\leq x_{qK}$ the points of discontinuity of $g^{(q)}$. They divide the interval $[0,1)$ into subintervals $[x_j,x_{j+1})$, $j=1,...,qK$ $(x_{qK+1}=x_1)$. It may happen that some of these intervals are empty. However, if $[x_j,x_{j+1})$ is not degenerated then $g^{(q)}|_{(x_j,x_{j+1})}$ is absolutely continuous and $(g^{(q)})'(x)=qS$ for all $x\in(x_j,x_{j+1})$. So we have
$$\int_{x_j}^{x_{j+1}}e^{2\pi ir(g^{(q)}+h^{(q)})}d\lambda=\int_{x_j}^{x_{j+1}}\frac{e^{2\pi ir h^{(q)}}}{2\pi ikg'^{(q)}}d(e^{2\pi ir g^{(q)}}).$$ 
We obtain by integrating by parts on each interval $(x_j,x_{j+1})$
\beab I_{r,q}=\frac{1}{2\pi ir}\sum_{j=1}^{Kq}\left (\frac{e^{2\pi ir(g_-^{(q)}(x_{j+1})+h_-^{(q)}(x_{j+1}))}}{g_-^{(q)'}(x_{j+1})} -\frac{e^{2\pi i(rg_+^{(q)}(x_j)+rh_+^{(q)}(x_j))}}{g_+^{(q)'}(x_j)}\right )\\
-\frac{1}{2\pi ir}\int_0^1e^{2\pi irg^{(q)}}d(e^{2\pi irh^{(q)}}/g'^{(q)}).\eeab
It follows that (cf. the proof of Theorem 3 from \cite{Lem-Pa})
$$|I_{r,q}|\leq\frac{1}{2\pi |r|}\left (\frac{2Kq}{Sq}+Var(e^{2\pi ikh^{(q)}}/g'^{(q)})\right )$$
and
\beab
Var(e^{2\pi irh^{(q)}}/g'^{(q)})\leq \sup\left (\frac{1}{|g'^{(q)}|}\right )Var(e^{2\pi irh^{(q)}})+Var(1/g'^{(q)})\\
\leq\frac{1}{Sq}2\pi|r|Var(h^{(q)})+\frac{1}{S^2q^2}Var(g'^{(q)})\\
\leq \frac{2\pi |r|}{S}Var(h)+\frac{1}{S^2q}Var(g').
\eeab
Finally
$$|I_{r,q}|\leq \frac{K}{\pi |r|S}+\frac{Var(h)}{S}+\frac{Var(g')}{2\pi |r|S^2q}=\frac{Var(h)}{S}+\mbox{O}(1/r)$$
uniformly in $q$. Since $Var(h)/S<c$, we get $\limsup_{q\to+\infty}|I_{r,q}|\leq c<1$ for $r$ large enough. \hfill $\square$

 \label{ratners} \section{Ratner's property for roof functions from a subfamily of $\mathcal{U}$}
 We will prove weak Ratner's property of the special flow $(T^f_t)_{t\in\R}$ with $f$ belonging to some subfamily of $\mathcal{U}$ (Section \ref{sec}) which will be specified below. From now on, we assume that $\alpha$ has bounded partial qoutients and let $C:=\sup\{a_n\}_{n=1}^{+\infty}+1$. First, let us recall the following result.
\begin{lemma}\label{abs} {\em(\cite{Fr-Lem}, Lemma 6.1.)}
 Let $T:\mathbb{T}\to \mathbb{T}$ be the rotation by $\alpha$ with bounded partial quotients and let $f:\mathbb{T}\to \R$ be absolutely continuous with zero mean. Then 
$$\sup_{0\leq n<q_{s+1}}\sup_{\|y-x\|<\frac{1}{q_s}}|f^{(n)}(y)-f^{(n)}x|\to 0,\;\; \mbox{as}\;s\to +\infty.$$
\end{lemma}
Let $f:\mathbb{T}\to \R$ be a positive function in $\mathcal{U}$ with the set of discontinuities $\{\beta_i\}_{i=1}^{+\infty}$ and the corresponding set of jumps $\{d_i\}_{i=1}^{+\infty}$. We make an additional assumption on the jumps of $f$. Namely, we assume that there exist $j\in \N$ and $\theta>0$ such that
\begin{equation}\label{theta}\sum_{i>j}|d_i|\leq\frac{|S|}{(2+\theta)(2C+1)((2C+1)^j+1)}.\end{equation}
\begin{remark} \label{speed} {\em There is some natural well-ordering on the set of points $\{\beta_i\}_{i=1}^{+\infty}$ optimal from the point of view of (\ref{theta}). We may assume that $|d_i|\geq|d_{i+1}|$ for $i=1,2,...$. Indeed, let $\sigma:\N\to\N$ be a permutation such that $|d_{\sigma(i)}|\geq |d_{\sigma(i+1)}|$ for all $i\in\N$. Let $S'=\sum_{i=1}^{+\infty}|d_i|$. Then $\sum_{i>j}|d_{\sigma(i)}|=S'-\sum_{i=1}^{j}|d_{\sigma(i)}|\leq S'-\sum_{i=1}^{j}|d_i|=\sum_{i>j}|d_i|<\frac{|S|}{(2+\theta)((2C+1)^j+1)}$. It follows that if (\ref{theta}) is satisfied for some permutation of jumps, it is also satisfied for the monotonic one.}
\end{remark}
Let us now define 
$$D:=\left\{\sum_{i=1}^{+\infty}m_id_i:\;0\leq m_i<2C+1 \right\}.$$ 
\begin{lemma} Under the above assumptions there exist $p,\eta>0$ such that \\
$(p-\eta,p+\eta)\subset(0,|S|)\setminus(D\cup-D)$, and consequently the set $(\sgn S)p-D$ is bounded away from 0.
\end{lemma} 
Proof: Let $j,\theta$ satisfy (\ref{theta}). We consider the set of sums $A:=\{\sum_{i=1}^jm_id_i$, $0\leq m_i<2C+1$, $i=1,...,j$\} and $B:=A\cup-A$ ($|B|\leq 2(2C+1)^j$) and let $\xi:=\frac{|S|}{(2+\theta)(2C+1)((2C+1)^j+1)}$. Then there exist $p\in(0,|S|)$ and $\eta>0$ such that $(p-\eta,p+\eta)\cap (B+(-\xi,\xi))=\emptyset$. Indeed, the number of points from $B$ in the interval $(0,|S|)$ is at most $(2C+1)^j$ (because $B$ is symmetric), so there exists an interval $I=[a,b]\subset (0,|S|)$ of length at least $\frac{|S|}{(2C+1)^j+1}$ and $B\cap I=\emptyset$. Hence, if we take $p$ to be midpoint of $I$, and any $0<\eta< \frac{|S|\theta}{2(2+\theta)((2C+1)^j+1)}$ then $p-\eta>a+\xi$, $p+\eta<b-\xi$, so $(p-\eta,p+\eta)\subset (a+\xi,b-\xi)$. But from the definition of $I$ we get $(a+\xi,b-\xi)\cap (B+(-\xi,\xi))=\emptyset$, and from (\ref{theta}), $D\cup-D\subset B+(-\xi,\xi)$. Finally, $(p-\eta,p+\eta)\subset(0,|S|)\setminus(D\cup-D)$. \hfill $\square$

\vspace{2ex}

Let us now define $P:=\cl(\sgn(S)p-D)$, then $0\notin P$ and $P$ is compact. The
proof of the following result contains some ideas from \cite{Fr-Lem}.

\begin{theorem}\label{ratner}
Let $T:\mathbb{T}\to\mathbb{T}$ be an irrational rotation by $\alpha$ with bounded partial quotients and $f:\mathbb{T}\to \R_+$ be bounded away from zero function in $\mathcal{U}$, satisfying (\ref{theta}). Then the special flow $(T_t^f)_{t\in\R}$ has the $R(\gamma,P)$ property for every $\gamma>0$. 
\end{theorem}
Proof: By (\ref{fac}),  $$f(x)=f_{pl}(x)+f_{ac}(x)=\sum_{i=1}^{+\infty}d_i\{x-\beta_i\}+c+f_{ac}(x),$$
for some constant $c\in \R$.\\
\indent Let us fix $\epsilon>0$ and $N\in \N$. Because the series $\sum_{i=1}^{+\infty} |d_i|$ converges, there exists a number $m(\epsilon)\in \N$ such that $$\sum_{i>m(\epsilon)}|d_i|<\frac{\epsilon}{4(2C+1)}.$$
Let $\kappa(\epsilon):=\frac{1}{m(\epsilon)(2C+1)}\min\{\epsilon/2pC,1/C^2\}$.\\
From Lemma \ref{abs} it follows that there exists $s_0$ such that for all $s\geq s_0$ 
\begin{equation}\label{ac}\sup_{0\leq n<q_{s+1}}\sup_{\|y-x\|<\frac{1}{q_s}}|f_{ac}^{(n)}(y)-f_{ac}^{(n)}(x)|<\epsilon/4,\end{equation}
and $\min\{\kappa(\epsilon),1\}q_{s_0}>N$. Let $\delta:=(\delta(\epsilon,N))=\frac{p}{|S|q_{s_0+1}}$ and take $x,y\in \T$ with $0<\|x-y\|<\delta$. Let $s$ denote the unique natural number such that 
\begin{equation}\label{xy}\frac{p}{|S|q_{s+1}}<\|x-y\|\leq \frac{p}{|S|q_s},\end{equation}
then $s\geq s_0$. Without loss of generality assume that $x<y$ and $S>0$, the proof for $S<0$ goes analogously. We will consider the sequence $\left(f_{pl}^{(n)}(y)-f_{pl}^{(n)}(x)\right)_{n\in\N}$. We have
\begin{multline} f_{pl}^{(n+1)}(y)-f_{pl}^{(n+1)}(x)=f_{pl}^{(n)}(y)-f_{pl}^{(n)}(x)+
\sum_{i=0}^{+\infty}d_i(\{y+n\alpha-\beta_i\}-\{x+n\alpha-\beta_i\})=\\
=f_{pl}^{(n)}(y)-f_{pl}^{(n)}(x)+\sum_{i=0}^{+\infty}d_i(y-x-\chi_{(x,y]}(\{\beta_i-n\alpha\})).\end{multline}
So setting $m_i=m_i(n)=\sum_{j=0}^{n-1} \chi_{(x,y]}(\{\beta_i-j\alpha\})$ we get that 
\begin{equation}\label{xy2}f_{pl}^{(n)}(y)-f_{pl}^{(n)}(x)=nS(y-x)-\bar{d_n},\end{equation}
where $$\overline{d_n}:=\overline{d_n}(x,y)=\sum_{i=1}^{+\infty}m_id_i.$$
Let us fix now $i\in\N$ and assume that $\{\beta_i-\ell\alpha\},\{\beta_i-r\alpha\}\in(x,y]$, with $0\leq \ell, r<q_{s+1}$, $\ell\neq r$. Then 
$$\|\{\beta_i-\ell\alpha\}-\{\beta_i-r\alpha\}\|\geq\|q_s\alpha\|>\frac{1}{2q_{s+1}}\geq\frac{1}{2Cq_s}.$$
So the number $m_i=m_i(q_{s+1})$ of discontinuities of the function $f_{pl}^{(q_{s+1})}$ in the interval $(x,y]$ which are of the form $\{\beta_i-j\alpha\}$ for some $0\leq j<q_{s+1}$, is smaller than 
$$2Cq_s|y-x|+1\leq 2C\frac{p}{|S|}+1\leq 2C+1.$$
Hence $(\overline{d_n})\in D$, $n=1,...,q_{s+1}$. Because of (\ref{xy2}) and (\ref{xy}), we have
$$f_{pl}^{(q_s)}(y)-f_{pl}^{(q_s)}(x)+\overline{d_{q_s}}=q_sS(y-x)\leq p$$
$$f_{pl}^{(q_{s+1})}(y)-f_{pl}^{(q_{s+1})}(x)+\overline{d_{q_{s+1}}}=q_{s+1}S(y-x)> p.$$
Moreover, for every natural $n$
$$0<f_{pl}^{(n+1)}(y)-f_{pl}^{(n+1)}(x)+\overline{d_{n+1}}-(f_{pl}^{(n)}(y)-f_{pl}^{(n)}(x)+\overline{d_{n}})=S(y-x)\leq \frac{p}{q_s}.$$
Hence, there exists an integer interval $I\subset[q_s,q_{s+1}]$ such that
$$\left|f_{pl}^{(n)}(y)-f_{pl}^{(n)}(x)+\overline{d_{n}}-p\right|<\frac{\epsilon}{2}\;\hbox{ for}\; n\in I$$
and
$$|I|\geq \min\left(\frac{\epsilon}{2p}q_s,q_{s+1}-q_s\right)\geq \min\left(\frac{\epsilon}{2pC},\frac{1}{C^2}\right)q_{s+1}.$$
Since $s\geq s_0$, by (\ref{ac}) we have
$$|f^{(n)}(y)-f^{(n)}(x)+\bar{d_{n}}-p|<\frac{3\epsilon}{4}\;\hbox{ for}\; n\in I.$$ 
Let us consider all $i$, $1\leq i\leq m(\epsilon)$. Then the number of discontinuities of $f^{(q_{s+1})}$ comming from $\beta_i-j\alpha$, $0\leq j<q_{s+1}$ in the interval $(x,y]$ is at most $(2C+1)m(\epsilon)$. 
It follows that we can split the interval $I$ into at most $m(\epsilon)(2C+1)$ integer intervals on which the discontinuities come only from $\beta_i-j\alpha$, with $i>m(\epsilon)$ and $j=0,...,q_{s+1}$. It follows that if $a,b$ are in any such integer interval then $m_i(a)=m_i(b)$, $i=1,...,m(\epsilon)$. Thus we can choose an integer subinterval $J\subset I$ such that if $\min J:=M$ and if $M,M+r\in J$ then by (\ref{theta})
\beab\label{nr}\overline{d}_{M+r}-\overline{d_M}=\sum_{i=1}^{+\infty}(m_i(M+r)-m_i(M))d_i\leq\sum_{i>m(\epsilon)} (2C+1)|d_i|\leq\frac{\epsilon}{4}
\eeab
and moreover
$$|J|>\frac{1}{m(\epsilon)(2C+1)}\min\left(\frac{\epsilon}{2pC},\frac{1}{C^2}\right)q_{s+1}=\kappa(\epsilon)q_{s+1}.$$
Set $d=d_M$. Then, in view of (\ref{nr})
$$|f^{(m)}(y)-f^{(m)}(x)-p-d|\leq |f^{(m)}(y)-f^{(m)}(x)-p-d_m|+|d_m-d|<\epsilon $$ $m\in J$.
Let now $J=[M,M+L]\cap \Z$, then
$$\frac{L}{M}\geq\frac{|J|}{q_{s+1}}\geq \kappa(\epsilon),\;\; M\geq q_s\geq q_{s_0}>N,\;\;L\geq |J|\geq \kappa(\epsilon)q_{s+1}>\kappa(\epsilon)q_{s_0}>N.$$
Because the special flow $\mathcal{T}^f$ is weakly mixing by Proposition \ref{weak mixing}, the automorphism $T_\gamma^f$ is ergodic for all $\gamma\neq 0$, and an application of Proposition \ref{condi} completes the proof. \hfill $\square$ 
\begin{remark}\label{stab} {\em The above theorem yields a stability result for the roof functions belonging to the set $\mathcal{U}$ and satisfying (\ref{theta}). Namely, consider a function $f\in \mathcal{U}$ with the set of discontinuity points $\{\beta_i\}_{i=1}^{+\infty}$ and the corresponding set of jumps $\{d_i\}_{i=1}^{+\infty}$ satisfying (\ref{theta}) with $\theta_f,j_f$ and a function $g\in BV_{a+j}(\T)$ such that 
\begin{equation}\label{perag} Var\;g\leq \min\left\{\frac{|S(f)|}{(2+\eta_g)(2C+1)((2C+1)^{j_f}+1)}, d_{j_f} \right\}\end{equation} for some $\eta_g>\frac{\theta_f+7}{\theta_f}$. Then $f+g\in \mathcal{U}$ and $f+g$ satisfies (\ref{theta}) for some $\theta_{f+g}, j_{f+g}$. Indeed notice that $|S(f+g)|=|S(f)+S(g)|\geq |S(f)|-|S(g)|\geq |S(f)|-Var\;g> 0$, so $f+g\in \mathcal{U}$. Moreover, let $\{r_i\}_{i=1}^{+\infty}$ denote the set of discontinuities of $f+g$ and $\{m_i\}_{i=1}^{+\infty}$ the set of jumps of the function $g$ and we reorder the set of discontinuities of $f+g$ to get a decreasing one (see Remark \ref{speed}). By setting $j_{f+g}=j_f$, taking any $0<\theta_{f+g}<\frac{1}{4+\theta_f+\eta_g}$ and using consecutively (\ref{perag}), (\ref{theta})+(\ref{perag}), the definition of $\theta_{f+g}$, and again (\ref{perag}), we obtain 
$$\sum_{i>j_{f+g}}|r_i|\leq\sum_{i>j_f}^{+\infty}|d_i|+\sum_{i=1}^{+\infty}|m_i|\leq$$
$$\frac{|S(f)|}{(2+\theta_f)(2C+1)((2C+1)^{j_f}+1)}+\frac{|S(f)|}{(2+\eta_g)(2C+1)((2C+1)^{j_f}+1)}\leq$$
$$\frac{|S(f)|}{(2+\theta_{f+g})(2C+1)((2C+1)^{j_f}+1)}-\frac{|S(f)|}{(2+\theta_{f+g})(2+\eta_g)(2C+1)((2C+1)^{j_f}+1)}\leq$$
$$\frac{|S(f)|}{(2+\theta_{f+g})(2C+1)((2C+1)^{j_f}+1)}-\frac{|S(g)|}{(2+\theta_{f+g})(2C+1)((2C+1)^{j_f}+1)}\leq$$
$$\frac{|S(f+g)|}{(2+\theta_{f+g})(2C+1)((2C+1)^{j_{f+g}}+1)}.$$ 
}
\end{remark}

\section{Absence of partial rigidity for the roof functions from a subfamily of $\mathcal{U}$}
In this section, we will show the absence of partial rigidity of special flows over irrational rotation by $\alpha$ having bounded partial qoutients ($C:=\sup\{a_n\}+1$) and under roof functions $f:\T\to \R_+, f\in \mathcal{U}$. We will use the following general lemma.
\begin{lemma}\label{part} {\em(\cite{Fr-Lem}, Lemma 7.1.)}
Let $S:(X,\mathscr{B},\mu)\to (X,\mathscr{B},\mu)$ be an ergodic automorphism and $g\in L^1(X,\mu)$ such that $g\geq c>0$. Suppose that the special flow $(S_t^g)_{t\in\R}$ is partially rigid along a sequence $(t_n)$, $t_n\to +\infty$. Then there exists $0<u\leq 1$ such that for every $0<\epsilon<c$ we have 
$$\liminf_{n\to+\infty}\mu\left(\{x\in X:\exists_{j\in \N} |g^{(j)}(x)-t_n|<\epsilon \}\right)\geq u.$$
\end{lemma}
Recall that  $\{\beta_i\}_{i=1}^{+\infty}$ denotes the set of discontinuity points and the corresponding set of jumps is $\{d_i\}_{i=1}^{+\infty}$. Let $C_1,C_2$ be as in Lemma \ref{dist}. Assume that there exists a function $\eta:\R_+\to\R_+$ such that $\liminf_{\epsilon\to 0}\eta(\epsilon)\epsilon=0$ and
\begin{equation}\label{eta}\sum_{i>\eta(\epsilon)}|d_i|<\frac{\epsilon}{\frac{C_2}{C_1}+1}.\end{equation}
We note that if such an $\eta$ exists it also exists for a monotonic permutation of the set of discontinuities of $f$ (cf. Remark \ref{speed}). Under the above assumption we prove the following result.

\begin{theorem}\label{rig} Assume $T:\T\to\T$ is an ergodic rotation by $\alpha$ having bounded partial quotients. Suppose $f\in\mathcal{U}$ is positive, bounded away from zero function satisfying (\ref{eta}). Then the special flow $(T_t^f)_{t\in\R}$ is not partially rigid.
\end{theorem}
Proof: Let $m,M$ be positive numbers such that $0<m\leq f(x)\leq M$ for every $x\in \T$. We will proceed by contradiction assuming that $(t_n)$ ($t_n\to +\infty$) is a partial rigidity time for $(T_t^f)_{t\in\R}$. By Lemma \ref{part}, there exists $0<u\leq 1$ such that for every $0<\epsilon<m$ we have 
\begin{equation}\label{pr}\liminf_{n\to+\infty}\lambda\left(\{x\in \T:\exists_{j\in \N} |f^{(j)}(x)-t_n|<\epsilon \}\right)\geq u.\end{equation}
 Let us fix $\epsilon>0$ such that  $\epsilon<\frac{m}{10}$ and 
\begin{equation}\label{aet}0<\eta(\epsilon)\epsilon<\frac{|S|m^2}{48M(m+Var\;f)+|S|m^2}u.\end{equation}
Since $f'\in L^1(\T,\lambda)$, there exists $0<\delta<\epsilon$ such that 

\begin{equation}\label{del} \lambda(A)<\delta \; \hbox{implies}\; \int_A|f'|d\lambda<\epsilon.\end{equation} Moreover by the Egorov theorem and the ergodicity of $T$ (and recalling that $S=\int_\T f'd\lambda$) it follows that there exist $A_\epsilon\subset \T$
with $\lambda(A_\epsilon)>1-\delta$ and $m_0\in\N$ such that 
\begin{equation}\label{mat}\begin{matrix}\frac{S}{2}\leq \frac{1}{k}f'^{(k)}(x),\;\hbox{if}\; S>0,\\
\\
\frac{S}{2}\geq \frac{1}{k}f'^{(k)}(x),\;\hbox{if}\; S<0
\end{matrix}\end{equation}
for all $k\geq m_0$ and $x\in A_\epsilon$.\\
\indent Take any $n\in \N$ such that $\frac{t_n}{2M}\geq m_0$, $t_n>2\epsilon$. Let us consider the set $J_{n,\epsilon}$ of all $j\in\N$ such that $|f^{(j)}(x)-t_n|<\epsilon$ for some $x\in\T$. For such $j$ and $x$ we have 
$$t_n+\epsilon\geq f^{(j)}(x)\geq mj\;\hbox{and}\; t_n-\epsilon<f^{(j)}(x)\leq Mj$$
 whence
$$\frac{t_n}{2M}\leq \frac{t_n-\epsilon}{M}<j<\frac{t_n+\epsilon}{m}\leq\frac{2t_n}{m}$$
for any $j\in J_{n,\epsilon}$; in particular, $j\in J_{n,\epsilon}$ implies $j\geq m_0.$\\
\indent Let now $j_n:=\max J_{n,\epsilon}$. Let us consider the points of discontinuity of $f^{(j_n)}$, i.e. $\{\beta_i-j\alpha\}$, $i=1,...,+\infty$, $0\leq j<j_n$. Consider first the points for which $i\leq\eta(\epsilon)$. They divide $\T$ into subintervals $I_1^n,...,I^n_{{\eta}(\epsilon)j_n}$. Some of these intervals can be empty. Note that the only discontinuities of $f^{(j_n)}$ which are contained in the interiors of $I_1^n,...,I^n_{{\eta}(\epsilon)j_n}$ come from the set $\{\beta_\ell\}_{\ell>\eta(\epsilon)}$.
 For every fixed $k\in\N$, the number of points of the form $\{\beta_k-r\alpha\}$, $0\leq r<j_n$ in the interval $I_i^n$ is not bigger than $\frac{C_2}{C_1}+1$. Indeed, if $\{\beta_k-t\alpha\},\{\beta_k-r\alpha\}\in I_i^n$, $t\neq r$ then (see Lemma \ref{dist})
$$\|\{\beta_k-t\alpha\}-\{\beta_k-r\alpha\}\|\geq \frac{C_1}{j_n},$$
so the number of such points is not bigger than $\frac{j_n}{C_1}|I_i^n|+1\leq \frac{C_2}{C_1}+1$.
Let us fix now $1\leq i\leq\eta(\epsilon)j_n$. For every $j\in J_{n,\epsilon}$ let $I^n_{i,j}$ stand for the minimal closed interval of $\overline{I_i^n}$ which includes the set $\{x\in I_i^n\;:\; |f^{(j)}(x)-t_n|<\epsilon\}$. Of course, $I^n_{i,j}$ may be empty.\\

\indent Moreover for every $k\geq 1$ let $m_k=m_k(n,w,i,j)=\sum_{s=0}^{w-1} \chi_{I^n_{i,j}}(\{\beta_k-s\alpha\})$. It follows that $m_k(n,j_n,i,j)\leq \frac{C_2}{C_1}+1$.\\
\indent If $I^n_{i,j}=[z_1,z_2]$ is not empty then using (\ref{inpa}) and (\ref{eta})
\begin{multline}\label{fj} \frac{1}{j}\left|\int_{I^n_{i,j}}f'^{(j)}d\lambda\right|\leq
\frac{|f^{(j)}(z_1^+)-f^{(j)}(z_2^-)|}{j}+\frac{|\sum_{k>\eta(\epsilon)}m_k(n,j,i,j)d_k|}{j}\leq \\
\frac{|f^{(j)}(z_1^+)-t_n|+|t_n-f^{(j)}(z_2^-)|}{j}+\frac{\sum_{k>\eta(\epsilon)}m_k(n,j,i,j)|d_k|}{j}\leq\\
\frac{2\epsilon}{j}+\frac{(\frac{C_2}{C_1}+1)\frac{\epsilon}{\frac{C_2}{C_1}+1}}{j}\leq \frac{3\epsilon}{j}\leq \frac{6M\epsilon}{t_n}
\end{multline}
because the only discontinuity points of $f^{(j)}$ in $I^n_{i,j}$ come from $\beta_k-s\alpha$, $k>\eta(\epsilon)$, $0\leq s<j$.\\
\indent Now, suppose that $x$ is the endpoint of $I^n_{i,j}$ and $y$ is the endpoint of $I^n_{i,j'}$ with $j\neq j'$. Then, by (\ref{inpa}) and (\ref{eta}), it follows that 
\begin{multline}\label{inte}\int_x^y|f'|^{(j_n)}d\lambda\geq \left|\int_x^yf'^{(j)}d\lambda\right|\geq |f^{(j)}(y^-)-f^{(j)}(x^+)|-
\sum_{k>\eta(\epsilon)}m_k(n,j_n,i,j)|d_k|\geq\\ |f^{(j)}(y^-)-f^{(j')}(y^-)|-|f^{(j')}(y^-)-t_n|-|f^{(j)}(x^+)-t_n|-\epsilon\geq\\
m-|f^{(j')}(y)-t_n|-|f^{(j)}(x)-t_n|-|f^{(j)}(x^+)-f^{(j)}(x)|-\epsilon\geq\\ m-4\epsilon\geq \frac{m}{2}.\end{multline}
  Let $K_i=\{j\in J_{n,\epsilon}: I^n_{i,j}\neq\emptyset\}$ and denote $|K_i|=r\geq 1$. Then there exist $r-1$ pairwise disjoint subintervals $H_t\subset I_i^n$, $t=1,...,r-1$, which are disjoint from $I^n_{i,j}$, $j\in K_i$ and fill up the space between those intervals. In view of (\ref{inte}), we have \begin{equation}\label{ht}\int_{H_t}|f'|^{(j_n)}d\lambda\geq \frac{m}{2} \end{equation}
for $t=1,...,r-1$. Therefore in view of (\ref{fj}) and (\ref{ht})
\begin{multline}\label{soj}\left|\sum_{j\in K_i}\int_{I^n_{i,j}}\frac {f'^{(j)}}{j}d\lambda\right|\leq r\frac{6M\epsilon}{t_n}
=\frac{6M\epsilon}{t_n}+\frac{12M\epsilon}{mt_n}(r-1)\frac{m}{2}\\ \leq 
\frac{6M\epsilon}{t_n}+\frac{12M\epsilon}{mt_n}\sum_{t=1}^{r-1}\int_{H_t}|f'|^{(j_n)}d\lambda \leq \frac{6M\epsilon}{t_n}+\frac{12M\epsilon}{mt_n}\int_{I^n_i}|f'|^{(j_n)}d\lambda.
\end{multline}
Since $\lambda(A_\epsilon^c)<\delta$, and $\frac{1}{j}<\frac{2M}{t_n}$ we have 
\begin{multline}\label{suy}\sum_{i=1}^{\eta(\epsilon)j_n}\sum_{j\in K_i}\int_{I^n_{i,j}\cap A^c_\epsilon}\frac{f'^{(j)}}{j}d\lambda\leq \frac{2M}{t_n}\sum_{i=1}^{\eta(\epsilon)j_n}\sum_{j\in K_i}\int_{I^n_{i,j}\cap A^c_\epsilon}|f'|^{(j_n)}d\lambda\leq \frac{2M}{t_n}\int_{ A^c_\epsilon}|f'|^{(j_n)}d\lambda\leq \\
\frac{2M}{t_n}j_n\epsilon\leq \frac{4M}{m}\epsilon.\end{multline}
Note that (see (\ref{pr}))
$$B_n:=\{x\in\T:\; \exists_{j\in\N}|f^{(j)}(x)-t_n|<\epsilon\}\subset \bigcup_{i=1}^{\eta(\epsilon)j_n}\bigcup_{j\in K_i}I^n_{i,j}. $$
Now, we can conclude as in \cite{Fr-Lem}, namely,
using (\ref{mat}), (\ref{soj}), (\ref{suy}) we have

$$\frac{|S|}{2}\lambda(B_n\cap A_\epsilon)\leq \sum_{i=1}^{\eta(\epsilon)j_n}\sum_{j\in K_i}\int_{I^n_{i,j}\cap A_\epsilon}\frac{|f'^{(j)}|}{j}d\lambda\leq \left|\sum_{i=1}^{\eta(\epsilon)j_n}\sum_{j\in K_i}\int_{I^n}\frac{f'^{(j)}}{j}d\lambda\right|+$$ $$\sum_{i=1}^{\eta(\epsilon)j_n}\sum_{j\in K_i}\int_{I^n_{i,j}\cap A_\epsilon^c}\frac{|f'^{(j)}|}{j}d\lambda\leq \eta(\epsilon)j_n\frac{6M\epsilon}{t_n}+\frac{12M\epsilon}{t_nm}\int_\T|f'|^{(j_n)}d\lambda +\frac{4M}{m}\epsilon$$
$$\leq\frac{12\eta(\epsilon)M\epsilon}{m}+\frac{4M}{m}\epsilon+\frac{24M\epsilon}{m^2}\|f'\|_{L^1}\leq 
\frac{24\eta(\epsilon)M}{m^2}(m+Varf)\epsilon
$$
Finally, using (\ref{aet}) we obtain
$$\lambda(B_n)\leq \lambda(B_n\cap A_\epsilon)+\lambda(A^c_\epsilon)<\frac{48\eta(\epsilon)M}{|S|m^2}(m+Varf)\epsilon+\epsilon<u.$$
The yields contradiction with (\ref{pr}) which completes the proof. \hfill $\square$

\vspace{2ex}

From Theorem \ref{ratner}, Theorem \ref{rig} and Lemma \ref{mild mixing} we get the following result. 
\begin{corollary}
Suppose that $T:\T\to \T$ is the rotation by an irrational number $\alpha$
with bounded partial quotients and $f:\T\to\R$ is a positive, bounded away from zero function from the set $\mathcal{U}$ satisfying conditions (\ref{theta}), (\ref{eta}). Then the special flow $(T_t^f)_{t\in\R}$ is mildly mixing.
\end{corollary}

Institute of Math.\\
Polish Academy of Scienes,\\
Sniadeckich 8,\\
00-950 Warszawa, Poland\\
adkanigowski@gmail.com

\end{document}